\documentclass[10pt,twoside]{article}
\usepackage{amsfonts}
\usepackage{fancyhdr}
\usepackage{titlesec}
\usepackage{cite}
\usepackage{ifthen}
\usepackage{pifont}
\usepackage{stmaryrd}
\usepackage{setspace}
\usepackage{indentfirst}
\usepackage{amsmath,amssymb,amscd,bbm,amsthm,mathrsfs,dsfont}
\input amssym.def
\usepackage{amssymb}
\usepackage{amsmath}
\usepackage{mathtools}
\usepackage{mathrsfs}
\usepackage{cite}
\usepackage{amsthm}
\usepackage{enumerate}
\theoremstyle{definition}
\numberwithin{equation}{section}
\newtheorem{proposition}{Proposition}[section]
\newtheorem{remark}{Remark}[proposition]
\newtheorem*{proof*}{Proof}
\newtheorem{lemma}{Lemma}[section]

\titleformat{\section}{\centering\large\bfseries}{\S\arabic{section}}{1em}{}
\newboolean{first}
\setboolean{first}{true}

\begin{document}

\setlength\abovedisplayskip{2pt}
\setlength\abovedisplayshortskip{0pt}
\setlength\belowdisplayskip{2pt}
\setlength\belowdisplayshortskip{0pt}

\title{\bf Some deviation inequalities for sums of negatively associated random variables  \author{WenCong Zhang}} \maketitle
 \footnote{Received: \today.}
 \footnote{MR Subject Classification: 0211 .}
\begin{center}
\begin{minipage}{135mm}

{\bf \small Abstract}.\hskip 2mm {\small
Let $\{X_i,i\geq1\}$ be a sequence of negatively associated random variables, and let $\{X_i^\ast,i\geq 1\}$ be a sequence of independent random variables such that $X_i^\ast$ and $X_i$ have the same distribution for each $i$. Denote by $S_k=\sum_{i=1}^{k}X_i$ and $S_k^\ast=\sum_{i=1}^{k}X_i^\ast$ for $k\geq 1$. The well-known results of Shao \cite{Shao2000} sates that $\mathbb{E}f(S_n)\leq \mathbb{E}f(S_n^\ast)$ for any nondecreasing convex function. Using this very strong property, we obtain a large variety of deviation inequalities for $S_n$.}
\end{minipage}
\end{center}

\thispagestyle{fancyplain} \fancyhead{}
\fancyhead[L]{\textit{Appl. Math. J. Chinese Univ.}\\
2017, 32(*): ***-***} \fancyfoot{} \vskip 10mm
\section{Introduction}
The concept of negatively associated (NA) is first introduced by Joag-Dev and Proschan \cite{joag-dev1983}.
A sequence of random variables $\{X_i,1\leq i\leq n\}$ is said to be NA if for every pair of disjoint subsets $A_1$ and $A_2$ of $\{1,2,\cdots,n\}$,
\[Cov\{f_1(X_i,i\in A_1),f_2(X_j,j\in A_2)\}\leq 0,\]
whenever $f_1$ and $f_2$ are coordinatewise increasing  and the covariance exists. A sequence of random variables $\{X_i,i\geq 1\}$ is said to be NA if for any $n\geq 2$ (where $n$ is a natural number), the sequence $X_1,X_2,\cdots,X_n$ is NA.

Since NA have a lot of applications in multivariate statistical analysis, reliability theory and percolation theory, many mathematicians have strong interest in it and discuss some property of NA sequences. Matula \cite{MATULA1992209} obtained a Kolmogorov type of upper bound inequality; Su et al. \cite{Su1997} gave some moment inequalities for NA sequences; Next, Liu et al. \cite{Liu1998} presented some probability and moment inequalities; Zhang et al. \cite{Liu19981} proved some Fuk-Negeav's inequalities; Gan et al. \cite{LIU199999} established the H\'{a}jeck-R\`{e}nyi inequality for NA sequences; Wang \cite{Wang2000} gave some exponential inequalities and a strong law of large numbers for NA random variables, etc.

In this paper, we are interested in establishing some new deviation inequalities for sums of NA random variables. The proof of our theorems are based on the following technical result of Shao \cite{Shao2000}. Throughout this paper, let $\{X_i,i\geq 1\}$ be a sequence of NA random variables, and let $\{X_i^\ast,i\geq 1\}$ be a sequence of independent random variables such that $X_i$ and $X_i^\ast$ have the same distribution for each $i\geq 1$. Denote by
\begin{equation}
 S_0=0,\qquad S_k=\sum_{i=1}^{k}X_i\qquad  \textrm{and}    \qquad S_k^\ast=\sum_{i=1}^{k}X_i^\ast\qquad \text{for all}\ \ k\geq 1.
\end{equation}
Then
\begin{equation}\label{s1:eq2}
\mathbb{E}f(S_n)\leq \mathbb{E}f(S_n^\ast)
\end{equation}
for any convex function $f$ on $\mathbb{R}$, whenever the expectation on the right hand side of (\ref{s1:eq2}) exists. If $f$ is a nondecreasing convex function, then
\begin{equation}\label{s1:eq3}
\mathbb{E}f(\max_{1\leq k\leq n}S_k)\leq \mathbb{E}f(\max_{1\leq k\leq n}S_k^\ast),
\end{equation}
whenever the expectation on the right hand side of (\ref{s1:eq3}) exists. Inequality (\ref{s1:eq2}) bridges the relationship between the sums of NA random variables and the sums of independent random variables. This enable us to extend many classical inequalities (cf.\cite{Lin2005}), such as the inequalities of Bernstain, Hoeffding, and  Nagaev for sums of independent (dependent) random variables to sums of NA random variables. In these inequalities, the tail probabilities are obtained by choosing some proper positive nondecreasing convex functions $f$ such that
\begin{equation}\label{s1:eq4}
\mathbb{P}(\max_{1\leq k\leq n}S_k\geq x)\leq \mathbb{E}f(\max_{1\leq k\leq n}S_k)/f(x).
\end{equation}
Applying (\ref{s1:eq3}) to the last inequality, we have
\begin{equation}\label{s1:eq5}
\mathbb{P}(\max_{1\leq k\leq n}S_k\geq x)\leq \mathbb{E}f(\max_{1\leq k\leq n}S_k^\ast)/f(x).
\end{equation}
The right-hand side of the last inequality is dominated by many classical moment inequalities with various conditions.
Similarly, by (\ref{s1:eq2}), (\ref{s1:eq5}) also holds when $\max_{1\leq k\leq n}S_k$ is replaced by $S_n$.
\section{Deviation inequalities}
In this section, we give some new deviation inequalities for NA sequences via Shao's comparison theorem on moment inequalities. Assume that the NA sequence $\{X_i,i\geq 1\}$ is centered, that is $\mathbb{E}X_i=0$ for all $i$. First we introduce a lemma which has been proved by Shao \cite{Shao2000}.
\begin{lemma}
If $\{T_i,1\leq i\leq n\}$ is a non-negative supermartingale, then, for any $0<\alpha<1$,
\begin{equation}
\mathbb{E}\max_{1\leq i\leq n}T_i^\alpha\leq \frac{(\mathbb{E}T_1)^\alpha}{1-\alpha}.
\end{equation}
This Lemma is important in our proof.
\end{lemma}
\subsection{Fuk-Nagaev type inequalities}
When the NA sequences have finite  $p$-th moments ($p\geq 2$), the following proposition gives a Fuk-Nagaev type inequality.

\begin{proposition}\label{pr2.1}
Assume $B_n(y):=\sum_{i=1}^{n}\mathbb{E}[X_i^2I_{\{X_i\leq y\}}]<\infty$ for any $y>0$. Then, for any $x, y>0$,
\begin{equation}\label{s2:eq28}
\mathbb{P}(\max_{1\leq k\leq n}S_k>x)\leq (1-\alpha)^{-1}H_n\biggl(\frac{\alpha x}{y},\frac{\sqrt{B_n(y)}}{y}\biggr)+\mathbb{P}(\max_{1\leq k\leq n}X_k>y),
\end{equation}
where
\begin{equation}\label{s2:eq2}
H_n(x,v)=\left\{\biggl(\frac{v^2}{x+v^2}\biggr)^{x+v^2}\biggl(\frac{n}{n-x}\biggr)^{n-x}\right\}^{\frac{\alpha n}{n+v^2}}I_{\{x\leq n\}}
\end{equation}
with the convention that $(+\infty)^0=1$ (which applies when $x=n$).
\end{proposition}
\begin{remark}
  According to Remark 2.1 of Fan et al.\cite{FAN20123545}, for any $x\geq 0$ and any $v>0$, it holds
\begin{eqnarray}\label{s2:eq3}\label{s2:eq4}
H_n(x,v)&\leq& B(x,v):=\biggl(\frac{v^2}{x+v^2}\biggr)^{\alpha (x+v^2)}e^{\alpha x}\\
\label{s2:eq4}&\leq &B_1(x,v):=\exp\left\{-\frac{\alpha x^2}{2(v^2+\frac{1}{3}x)}\right\}.
\end{eqnarray}
Note that (\ref{s2:eq3}) and (\ref{s2:eq4}) are respectively known as Bennett's and Bernstein's bounds. Then, inequality (\ref{s2:eq28}) also implies the inequalities of  Bennett and Bernstein: for any $x,y>0$,
\begin{eqnarray}
\mathbb{P}(\max_{1\leq k\leq n}S_k>x)&\leq& (1-\alpha)^{-1}B\biggl(\frac{\alpha x}{y},\frac{\sqrt{B_n(y)}}{y}\biggr)+\sum_{i=1}^{n}\mathbb{P}(X_i>y)\\
&\leq& (1-\alpha)^{-1}B_1\biggl(\frac{\alpha x}{y},\frac{\sqrt{B_n(y)}}{y}\biggr)+\sum_{i=1}^{n}\mathbb{P}(X_i>y).
\end{eqnarray}
For similar results, we refer to Su et al. \cite{Su1997}, Shao \cite{Shao2000} and Zhang and Liu \cite{Liu19981}.
\end{remark}
\begin{proof*}
For given $y>0$, let
\begin{eqnarray*}
&Y_i=\min\{X_i,y\},\qquad\qquad \tilde{S_k}=\sum_{i=1}^{k}Y_i,&\\
&\tilde{S_k^\ast}=\sum_{i=1}^{k}Y_i^\ast,\qquad \textrm{for all}\quad i\geq 1.&
\end{eqnarray*}
Then $\{Y_i,1\leq i\leq n\}$ is also a NA sequence. It is obviously that
\begin{equation*}
\mathbb{P}(\max_{1\leq k\leq n}S_k\geq x)\leq \mathbb{P}(\max_{1\leq k\leq n}\tilde{S_k}\geq x)+\mathbb{P}(\max_{1\leq k\leq n}X_k>y).
\end{equation*}
Let
\[T_k=\exp\big\{\frac{t}{y}\tilde{S_k^\ast}-nf\big(\frac{t}{y},\frac{B_n(y)}{y^2n}\big)\big\},\]
where $t>0$. Using Lemma 3.3 and Lemma 3.2 in Fan et al.\cite{FAN20123545}, we have
\begin{equation*}
\begin{aligned}
\mathbb{E}e^{\frac{t}{y}\tilde{S}_k^\ast}&= \exp\Big\{\sum_{i=1}^{k}\ln\mathbb{E}\big[e^{\frac{t}{y}\tilde{Y}_i^\ast}\big]\Big\}\leq \exp\Big\{kf\big(\frac{t}{y},\frac{B_k^\ast(y)}{k}\big)\Big\}\\
&\leq \exp\Big\{kf\big(\frac{t}{y},\frac{B_n(y)}{n}\big)\Big\},
\end{aligned}
\end{equation*}
where \[T_k=\exp\big\{\frac{t}{y}\tilde{S}_k^\ast-kf\big(\frac{t}{y},\frac{B_n(y)}{n}\big)\big\}.\]
Hence, $\{T_k,k\geq 1\}$ is a supermartingale. Applying $f(x)=e^{t\alpha x}$ for any $0<\alpha<1$, to (\ref{s1:eq5}), we obtain
\begin{equation*}
\begin{aligned}
\mathbb{P}(\max_{1\leq k\leq n }\tilde{S_k}\geq x)&=\mathbb{P}(\frac{\max_{1\leq k\leq n }\tilde{S_k}}{y}\geq \frac{x}{y})
\leq \frac{\mathbb{E}e^{\frac{t\alpha}{y}\max_{1\leq k\leq n}S_k^\ast}}{e^{t\alpha x}}\\
&\leq e^{-t\alpha \frac{x}{y}}\mathbb{E}\Big\{\max_{1\leq k\leq n}T_ke^{nf\big(\frac{t}{y},\frac{B_n(y)}{y^2n}\big)}\Big\}^\alpha\\
&\leq \exp\Big\{-t\alpha \frac{x}{y}+\alpha nf\big(\frac{t}{y},\frac{B_n(y)}{y^2n}\big)\mathbb{E}\Big\{\max_{1\leq k\leq n}T_k\Big\}^\alpha\\
&\leq (1-\alpha)^{-1}\exp\Big\{-t\alpha x+\alpha nf\big(\frac{t}{y},\frac{B_n(y)}{y^2n}\big)\Big\}.
\end{aligned}
\end{equation*}
Taking
\[t=\frac{1}{1+\frac {B_n(y)} {n}}\ln\frac{1+\frac {\alpha x} {B_n(y)}}{1-\frac \alpha x n},\]
then
\begin{equation*}
\mathbb{P}(\max_{1\leq k\leq n }S_k\geq x)\leq (1-\alpha)^{-1}H_n(\frac{\alpha x}{y},\frac{\sqrt{B_n(y)}}{y}),
\end{equation*}
which gives (\ref{s2:eq28}). \hfill\qed
\end{proof*}

Next we consider the case where the random variables $\{X_i,i\geq 1\}$ have only a weak moment of order $p\geq 2$. Recall the weak moment of order $p$ is defined by
\begin{equation}\label{s2:eq35}
||Z||_{w,p}^p=\sup_{x>0}x^p\mathbb{P}(|Z|>x)
\end{equation}
for any real-valued random variable $Z$ and any $p\geq 1$. By Proposition \ref{pr2.1}, it follows that:
\begin{proposition}\label{pr2.2}
Let $p\geq 2$. Assume
\[B_n:=\sum_{i=1}^{n}\mathbb{E}[X_i^2]<\infty \ \ \ \    \textrm{and}  \ \ \ \   A(p):=\sum_{i=1}^{n}||X_i||_{w,p}^p<\infty.\]
Then, for any $x,y>0$,
\begin{equation}\label{s2:eq36}
\mathbb{P}(\max_{1\leq k\leq n}S_k>x)\leq (1-\alpha)^{-1}H_n\big(\frac{\alpha x}{y},\frac{\sqrt{B_n}}{y}\big)+\frac{A(p)}{y^p},
\end{equation}
where $H_n(x,v)$ is defined by (\ref{s2:eq2}).
\end{proposition}
\begin{remark}
Assume that the random variables $\{X_i,i\geq 1\}$ have a weak moment of order $p>2$. Note that $H_n(x,v)\leq B_1(x,v)$. Taking
\[y=\frac{3nx}{2p\ln n} \]
in inequality (\ref{s2:eq36}), we infer that, for any $x>0$,
\begin{equation}\label{s2:eq37}
\mathbb{P}(\max_{1\leq k\leq n}S_k>nx)\leq C_x\frac{(\ln n)^p}{n^{p-1}}
\end{equation}
for some positive $C_x$ not depending on $n$.
\end{remark}
If the weak $p$-th moments of the random variables $\{X_i,i\geq 1\}$ are strengthen to the $p$-th moments ($p\geq 2$), then we have the following Fuk-type inequality (cf. Corollary 3$'$ of Fuk \cite{Fuk1973}).
\begin{proposition}
Let $p\geq 2$. Assume
\begin{equation}\label{s2:eq38}
V_n:=\sum_{i=1}^{n}\mathbb{E}[|X_i|^p]<\infty.
\end{equation}
Then, for any $x>0$,
\begin{equation}\label{s2:eq39}
\mathbb{P}(\max_{1\leq k\leq n}S_k>x)\leq (1-\alpha)^{-1}\Big[\biggl(1+\frac{2}{p}\biggr)^p\frac{V_n}{\alpha^px^p}+\exp\biggl\{-\alpha\frac{2}{(p+2)^2e^p}\frac{x^2}{B_n}\biggr\}\Big],
\end{equation}
where $B_n$ is defined by Proposition \ref{pr2.2}.
In particular, we have
\begin{equation}\label{s2:eq102}
\mathbb{P}(\max_{1\leq k\leq n}S_k>x)\leq 2^{p+1}\biggl(1+\frac{2}{p}\biggr)^p\frac{V_n}{x^p}+2\exp\biggl\{-\frac{x^2}{(p+2)^2e^p B_n}\biggr\}.
\end{equation}
\end{proposition}
\begin{remark}Consider the case that $\{X_i,i\geq 1\}$ is a stationary sequence. Since $B_n$ and $V_n$ are of order $n$ as $n\rightarrow \infty$, it easy to see that the sub-Gaussian term
\[\exp\biggl\{-\frac{2\alpha }{(p+2)^2e^p}\frac{x^2}{B_n}\biggr\}\]
is decreasing at an exponential order, and that the polynomial term
\[2\biggl(1+\frac{2}{p}\biggr)^p\frac{V_n}{(x\alpha n)^p}\]
is of order $n^{1-p}$. Thus, for any $x>0$ and all $n$,
\[\mathbb{P}(|S_n|>nx)\leq \frac{C_x}{n^{p-1}}\]
for some positive $C_x$ not depending on $n$. Thus the order $\frac{(\ln n)^p}{n^{p-1}}$ in (\ref{s2:eq37}) is refined to $\frac{1}{n^{p-1}}$.
\end{remark}
\begin{proof*}
For given $x,t,y_1,y_2,\cdots,y_n>0$, let $y\geq \max\{y_1,\cdots,y_n\}$ and
\[
\tilde{X_i}=\min{\{X_i,y_i\}},\qquad\qquad \tilde{S_k}=\sum_{i=1}^{k}\tilde{X_i},\]
\[
\tilde{S_k^\ast}=\sum_{i=1}^{k}{\tilde{X_i^\ast}}.\qquad k=1,2,\cdots,n.\]
Then $\{\tilde{X_i},i\geq 1\}$ is also a NA sequence. It is easy to see that
\begin{equation*}
\mathbb{P}(\max_{1\leq k\leq n}S_k\geq x)\leq \sum_{i=1}^{n}\mathbb{P}(X_i\geq y_i)+\mathbb{P}(\max_{1\leq k\leq n}\tilde{S_k}\geq x).
\end{equation*}
Let
\[T_k=\exp\big\{t\tilde{S_k^\ast}-\big(\frac{e^{ty}-1-ty}{y^p}\big)V_k-\frac12e^pB_kt^2\big\}.\]
By Lemma 2 in Fuk \cite{Fuk1973}, we have
\[\mathbb{E}e^{t\tilde{S_n^\ast}}\leq \exp\Big\{\Big(\frac{e^{ty}-1-ty}{y^p}\Big)V_n+\frac12e^pB_nt^2\Big\}.\]
Hence, $\{T_k,k\geq 1\}$ is a supermartingale. Applying $f(x)=e^{t\alpha x}$ for any $0<\alpha<1$, to (\ref{s1:eq5}), we deduce that
\begin{equation}\label{s2:eq87}
\begin{aligned}
\mathbb{P}(\max_{1\leq k\leq n}\tilde{S_k}\geq x)&\leq \frac{\mathbb{E}e^{t\alpha\max_{1\leq k\leq n}\tilde{S_k^\ast}}}{e^{t\alpha x}}\\
&\leq e^{-t\alpha x}\mathbb{E}\Big\{\max_{1\leq k\leq n}T_ke^{\Big(\frac{e^{ty}-1-ty}{y^p}\Big)V_k+\frac12e^pB_kk^2}\Big\}^\alpha\\
&\leq \exp\Big\{-t\alpha x+\alpha\Big(\Big(\frac{e^{ty}-1-ty}{y^p}\Big)V_n+\frac12e^pB_nt^2\Big)\Big\}\mathbb{E}\Big\{\max_{1\leq k\leq n}T_k\Big\}^\alpha\\
&\leq (1-\alpha)^{-1}\exp\Big\{-t\alpha x+\alpha\Big(\Big(\frac{e^{ty}-1-ty}{y^p}\Big)V_n+\frac12e^pB_nt^2\Big)\Big\}\\
&=(1-\alpha)^{-1}\exp\Big\{\alpha\big(f_1(t)+f_2(t)\big)\Big\},
\end{aligned}
\end{equation}
where
\[f_1(t)=\Big(\frac{e^{ty}-1-ty}{y^p}\Big)V_n-\mu tx\qquad \textrm{and} \qquad f_2(t)=\frac12e^pB_nt^2-\lambda tx.\]
Let
\[t_1=\max\Big\{\frac{p}{y},\frac{\ln(\frac{\mu xy^{p-1}}{V_n}+1)}{y}\Big\}\qquad \textrm{and}\qquad t_2=\frac{\lambda x}{e^pB_n}.\]
If
\[t_1\geq t_2,\]
then
\begin{equation*}
\begin{aligned}
\mathbb{P}(\max_{1\leq k\leq n}\tilde{S_k}\geq x)&\leq (1-\alpha)^{-1}\exp\Big\{\alpha\big(f_1(t_2)+f_2(t_2)\big)\Big\}\\
&\leq (1-\alpha)^{-1}\exp\Big\{-\frac{\alpha\lambda^2x^2}{2e^pB_n}\Big\}.
\end{aligned}
\end{equation*}
If
\[t_1<t_2,\]
then
\begin{equation*}
\begin{aligned}
\mathbb{P}(\max_{1\leq k\leq n}\tilde{S_k}\geq x)&\leq (1-\alpha)^{-1}\exp\Big\{\alpha\big(f_1(t_1)+f_2(t_1)\big)\Big\}\\
&\leq (1-\alpha)^{-1}\exp\Big\{-\alpha\mu\frac xy\ln\big(\frac{\mu xy^{p-1}}{V_n}+1\big)\Big\}.
\end{aligned}
\end{equation*}
So we have
\begin{equation*}
\mathbb{P}(\max_{1\leq k\leq n}\tilde{S_k}\geq x)\leq (1-\alpha)^{-1}\exp\Big\{-\alpha\mu\frac xy\ln\big(\frac{\mu xy^{p-1}}{V_n}+1\big)\Big\}+(1-\alpha)^{-1}\exp\Big\{-\frac{\alpha\lambda^2x^2}{2e^pB_n}\Big\}.
\end{equation*}
Setting $y=y_1=y_2=\cdots=y_n=\alpha\mu x$. Since
\[\sum_{i=1}^{n}\mathbb{P}(X_i\geq y_i)\leq \frac{\sum_{i=1}^{n}\mathbb{E}X_i^p}{y_1^p}\leq \frac{V_n}{\alpha^p\mu^px^p},\]
and
\[\exp\Big\{-\alpha\frac{\mu x}{y}\ln\big(\frac{\mu xy^{p-1}}{V_n}+1\big)\Big\}=\frac{V_n}{\mu^px^p\alpha^{p-1}+V_n}\leq \frac{V_n}{\mu^px^p\alpha^{p-1}}=\frac{\alpha V_n}{\mu^px^p\alpha^p},\]
where $\lambda=\frac{2}{p+2}$ and $\mu=1-\lambda$. Then we obtain (\ref{s2:eq39}). \hfill\qed
\end{proof*}
\subsection{Semi-exponential bound}
When the random variables $\{X_i,i\geq 1\}$ have semi-exponential moments, the following proposition holds. This proposition can be compared to the corresponding results in Fan et al.\cite{FAN2017538} for martingales.
\begin{proposition}\label{pr2.6}
Let $p\in(0,1)$. Assume that there exists a positive constant $K_n\geq 1$ such that
\begin{equation}\label{s2:eq58}
\sum_{i=1}^{n}\mathbb{E}[X_i^2\exp\{|X_i|^p\}]\leq K_n.
\end{equation}
Then, for any $x>0$,
\begin{eqnarray}\label{s2:eq17}
\mathbb{P}(\max_{1\leq k\leq n}S_k\geq x)&\leq& \begin{cases}
2(1-\alpha)^{-1}\exp\left\{-\frac{\alpha^2x^2}{2K_n}\right\}\qquad&\text{if}\quad 0\leq \alpha x\leq K_n^{\frac{1}{2-p}}\\
2(1-\alpha)^{-1}\exp\left\{-\frac{\alpha^px^p}{2}\right\}\qquad&\text{if}\quad  \alpha x\geq K_n^{\frac{1}{2-p}}
\end{cases}\\
&\leq& 2(1-\alpha)^{-1}\exp\left\{-\frac{\alpha^2x^2}{2(K_n+\alpha^{2-p}x^{2-p})}\right\}.
\end{eqnarray}
In particular, we have
\begin{eqnarray}\label{s2:eq20}
\mathbb{P}(\max_{1\leq k\leq n}S_k\geq x)&\leq& \begin{cases}
4\exp\left\{-\frac{x^2}{8K_n}\right\}\qquad&\text{if}\quad 0\leq x\leq 2K_n^{\frac{1}{2-p}}\\
4\exp\left\{-\frac{x^p}{2^{p+1}}\right\}\qquad&\text{if}\quad   x\geq 2K_n^{\frac{1}{2-p}}
\end{cases}\\
&\leq& 4\exp\left\{-\frac{x^2}{8K_n+2^{p+1}x^{2-p}}\right\}.
\end{eqnarray}
\end{proposition}
\begin{remark}
It is interesting to see that for moderate $\alpha x\in (0,K^{\frac{1}{2-p}})$, the bound (\ref{s2:eq17}) is a sub-Gaussian bound and is of order
\[
\exp\left\{-\frac{\alpha^2x^2}{2K_n}\right\}.
\]
For all $x\geq K_n^{\frac{1}{2-p}}$, bound (\ref{s2:eq17}) is a semi-exponential bound and is of order
\[
\exp\left\{-\frac{\alpha^px^p}{2}\right\}.
\]
In particular, in the stationary case, there exists a positive constant $c$ such that, for any $x>0$,
\begin{equation}\label{s2:eq60}
\mathbb{P}(\max_{1\leq k\leq n}S_k\geq nx)\leq 2(1-\alpha)^{-1}\exp\{-c\alpha^px^pn^p\},
\end{equation}
where the constant $c$ does not depend on $n$.
\end{remark}
\begin{proof*}
For given $x,y,t>0$, denote
\[
\tilde{X_i}=\min{\{X_i,y\}},\quad \tilde{S_k}=\sum_{i=1}^{k}\tilde{X_i},\ \ \
\tilde{S_k^\ast}=\sum_{i=1}^{k}{\tilde{X_i^\ast}},\quad k=1,2,\cdots,n.\]
It is easy to see that
\begin{equation*}
\begin{aligned}
\mathbb{P}(\max_{1\leq k\leq n}S_k\geq x)&\leq \mathbb{P}(\max_{1\leq k\leq n}X_k\geq y)+\mathbb{P}(\max_{1\leq k\leq n}\tilde{S_k}\geq x)\\
&\leq \sum_{i=1}^{n}\mathbb{P}(X_i\geq y)+\mathbb{P}(\max_{1\leq k\leq n}\tilde{S_k}\geq x).
\end{aligned}
\end{equation*}
Notice that $\{\tilde{X_i},i\geq 1\}$ is also a NA sequence. Let $\psi_n(t)=\sum_{i=1}^{n}\log\mathbb{E}e^{t \tilde{X_i^\ast}}$ and $t=y^{p-1}$. By Lemma 4.1 in Fan et al. \cite{FAN2017538} and the inequality $\log(1+t)\leq t$ for all $t\geq 0$, we obtain
\begin{equation*}
\begin{aligned}
\mathbb{E}e^{t\tilde{S_n^\ast}}&=\exp\{\psi_n(t)\}\\
&\leq \exp\Big\{\sum_{i=1}^{n}\log\Big(1+\frac{\lambda^2}{2}\mathbb{E}[{\tilde{X_i^\ast}}^2\exp\{\lambda y^{1-p}({\tilde{X_i^\ast}}^{+})^p\}]\Big)\Big\}\\
&\leq \exp\Big\{\sum_{i=1}^{n}\frac{\lambda^2}{2}\mathbb{E}[{\tilde{X_i^\ast}}^2\exp\{\lambda y^{1-p}({\tilde{X_i^\ast}}^{+})^p\}]\Big\}\\
&\leq \exp\Big\{\frac12 y^{2p-2}K_n\Big\}.
\end{aligned}
\end{equation*}
Let
\[T_k=\exp\big\{t\tilde{S_k^\ast}-\frac12 y^{2p-2}K_k\big\}.\]
Then $\{T_k,k\geq 1\}$ is a supermartingale. Applying $f(x)=e^{t\alpha x}$ for any $0<\alpha<1$, to (\ref{s1:eq5}), we have
\begin{equation*}
\begin{aligned}
\mathbb{P}(\max_{1\leq k\leq n}\tilde{S_k}\geq x)&\leq\frac{\mathbb{E}e^{t\alpha\max_{1\leq k\leq n}\tilde{S_k^\ast}}}{e^{t\alpha x}}\\
&\leq e^{-t\alpha x}\mathbb{E}\Big\{\max_{1\leq k\leq n}T_ke^{\frac12 y^{2p-2}K_k}\Big\}^\alpha\\
&\leq \exp\Big\{-t\alpha x+\alpha \frac12 y^{2p-2}K_n\Big\}\mathbb{E}\Big\{\max_{1\leq k\leq n}T_k\Big\}^\alpha\\
&\leq (1-\alpha)^{-1}\exp\Big\{-t\alpha x+\alpha \frac12 y^{2p-2}K_n\Big\}\\
&= (1-\alpha)^{-1}\exp\Big\{-y^{p-1}\alpha x+\alpha \frac12 y^{2p-2}K_n\Big\}
\end{aligned}
\end{equation*}
and
\begin{equation*}
\begin{aligned}
\sum_{i=1}^{n}\mathbb{P}(X_i\geq y)&\leq \frac{1}{y^2}\exp\{-y^p\}\sum_{i=1}^{n}\mathbb{E}[X_i^2\exp\{|X_i|^p\}]\\
&\leq \frac{K_n}{y^2}\exp\{-y^p\}\\
&< (1-\alpha)^{-1}\frac{K_n}{y^2}\exp\{-y^p\}.
\end{aligned}
\end{equation*}
Then
\begin{equation}\label{s2:eq21}
\mathbb{P}(\max_{1\leq k\leq n}S_k\geq x)\leq (1-\alpha)^{-1}\frac{K_n}{y^2}\exp\Big\{-y^p\Big\}+(1-\alpha)^{-1}\exp\Big\{-y^{p-1}\alpha x+\alpha \frac12 y^{2p-2}K_n\Big\}.
\end{equation}
Taking
\begin{equation*}
y=\begin{cases}
 (\frac{K_n}{\alpha x})^{\frac{1}{1-p}}\qquad&\textrm{if}\quad 0\leq \alpha x\leq K_n^{\frac{1}{2-p}}\\
 \alpha x\qquad& \text{if}\quad \alpha x\geq K_n^{\frac{1}{2-p}}
 \end{cases}
\end{equation*}
to (\ref{s2:eq21}), then we obtain the desired inequality. \hfill\qed
\end{proof*}
\subsection{Qualitative results when $\mathbb{E}[e^{a|X_i|^p}]<\infty$ for $p>1$}
Let $\{X_i,i\geq 1\}$ be a sequence of non-degenerate NA random variables in this subsection. This proposition can be compared with Liu and Watbled \cite{LIU20093101} for martingales.
\begin{proposition}\label{pr2.5}
Let $p>1$. Assume that there exists a constant $a>0$ such that
\begin{equation}\label{s2:eq82}
K:=\sum_{i=1}^{n}\mathbb{E}[\exp\{a|X_i|^p\}]<\infty.
\end{equation}
Let $q$ be the conjugate exponent of $p$ and let $\tau>0$ be such that
\[(q\tau)^{\frac{1}{q}}(pa)^{\frac{1}{p}}=1.\]
Then, for any $\tau_1>\tau$, there exist some positive numbers $t_1$, $x_1$, $A$ and $B$, depending only on $a$, $K$ and $p$, such that
\begin{equation}\label{s2:eq83}
\mathbb{E}[e^{tS_n}]\leq\begin{cases}
\exp\{n\tau_1t^q\}\qquad &\text{if}\quad t\geq t_1\\
\exp\{nAt^2\}\qquad&\text{if}\quad 0\leq t\leq t_1,
\end{cases}
\end{equation}
and for any $x>0$,
\begin{equation}\label{s2:eq84}
\mathbb{P}(\max_{1\leq k\leq n}S_k\geq x)\leq\begin{cases}
(1-\alpha)^{-1}\exp\{-a_1\frac{\alpha x^p}{n^{p-1}}\}\qquad &\text{if}\quad x\geq nx_1\\
(1-\alpha)^{-1}\exp\{-B\frac{\alpha x^2}{n}\}\qquad&\text{if}\quad 0\leq x\leq nx_1.
\end{cases}
\end{equation}
In particular, we have
\begin{equation}
\mathbb{P}(\max_{1\leq k\leq n}S_k\geq x)\leq\begin{cases}
2\exp\{-a_1\frac{x^p}{2n^{p-1}}\}\qquad &\text{if}\quad x\geq nx_1\\
2\exp\{-B\frac{x^2}{2n}\}\qquad&\text{if}\quad 0\leq x\leq nx_1,
\end{cases}
\end{equation}
where $a_1$ is such that $(q\tau_1)^{\frac{1}{q}}(pa_1)^{\frac{1}{p}}=1$.
\end{proposition}
\begin{remark}
Let us comment on Proposition \ref{pr2.5}.
\begin{itemize}
  \item [1.] Assume that (\ref{s2:eq84}) is satisfied for some $p>1$. From Proposition \ref{pr2.5}, we infer that for any $x>0$, one can find a positive constant $c_x$ not depending on $n$ such that
\begin{equation}\label{s2:eq85}
\mathbb{P}(\max_{1\leq k\leq n}S_k\geq nx)\leq (1-\alpha)^{-1}\exp\{-c_xn\}.
\end{equation}
Moreover, for $x$ large enough, one can take $c_x=a_1x^p$.
\item [2.] In particular, if $p=2$, we have the following sub-Gaussian bound. Assume that there exists a constant $a>0$ such that
\begin{equation}\label{s2:eq86}
K:=\sum_{i=1}^{n}\mathbb{E}[\textrm{exp}\{a|x_i|^2\}]< \infty.
\end{equation}
Since $1/p+1/q=1$, then $q=2$. So we have the following result,
\begin{equation}
\begin{aligned}
\mathbb{P}(\max_{1\leq k\leq n}S_k \geq nx)&\leq\begin{cases}
(1-\alpha)^{-1}\exp\{-a_1\alpha nx^2\} \qquad &\text{if}\quad x\geq {x_1} \\
(1-\alpha)^{-1}\exp\{-B\alpha nx^2\} \qquad &\text{if}\quad  0\leq  x\leq x_1.
\end{cases}\\
\end{aligned}
\end{equation}
\item [3.] By the inequality (\ref{s2:eq83}), we can get
\begin{equation}
\mathbb{P}(S_n\geq x)\leq\begin{cases}
\exp\{-a_1\frac{x^p}{n^{p-1}}\}\qquad &\text{if}\quad x\geq nx_1\\
\exp\{-B\frac{x^2}{n}\}\qquad&\text{if}\quad 0\leq x\leq nx_1,
\end{cases}
\end{equation}
\end{itemize}
\end{remark}
\begin{proof*}
By Lemmas 3.5 and 3.3 in Liu and Watbled \cite{LIU20093101}, we see that for $\alpha=K(\frac2a)^{\frac{1}{p-1}}$,
\[\mathbb{E}e^{t|X_i^\ast|}\leq 1+K+\alpha t^qe^{\tau t^q}.\]
Let $\tau_1>\tau$. Then there exists $t_1>\frac a2$ sufficiently large such that for any $t\geq t_1$,
\[\mathbb{E}e^{t|X_i^\ast|}\leq e^{\tau_1t^q}.\]
So, by (\ref{s1:eq2}), we have
\[\mathbb{E}e^{tS_n}\leq \mathbb{E}e^{tS_n^\ast}\leq \mathbb{E}e^{t|S_n^\ast|}\leq e^{n\tau_1t^q}.\]
Let
\[T_k=\exp\Big\{tS_k^\ast-k\tau_1t^q\Big\}.\]
Then $\{T_k,k\geq 1\}$ is a supermartingale. Applying $f(x)=e^{t\alpha x}$ for any $t>0$, to (\ref{s1:eq5}) and using Lemmas 2.3 and 3.4 in Liu and Watbled \cite{LIU20093101}, we obtain
\begin{equation*}
\begin{aligned}
\mathbb{P}(\max_{1\leq k\leq n}S_k\geq x)&\leq \frac{\mathbb{E}e^{t\alpha \max_{1\leq k\leq n}S_k^\ast}}{e^{t\alpha x}}\\
&\leq e^{-t\alpha x}\mathbb{E}\Big\{\max_{1\leq k\leq n}T_ke^{k\tau_1t^q}\Big\}^\alpha\\
&\leq \exp\Big\{-t\alpha x+\alpha n\tau_1t^q\Big\}\mathbb{E}\Big\{\max_{1\leq k\leq n}T_k\Big\}^\alpha\\
&\leq (1-\alpha)^{-1}\exp\Big\{-t\alpha x+\alpha n\tau_1t^q\Big\}\\
&\leq (1-\alpha)^{-1}\exp\Big\{-\alpha n\sup_{t}\big\{\frac{tx}{n}-\tau_1t^q\big\}\Big\}\\
&=(1-\alpha)^{-1}\exp\Big\{-a_1\frac{\alpha x^p}{n^{p-1}}\Big\}\qquad \textrm{if}\quad x\geq nx_1=q\tau_1t_1^{q-1}.
\end{aligned}
\end{equation*}
On the other hand, notice that $$\mathbb{E}[e^{a|X_i^\ast|}]\leq K_1:=e^a+K,$$ so by Theorem 2.1 in Liu and Watbled \cite{LIU20093101},
\[\mathbb{E}[e^{tS_n^\ast}]\leq \exp\Big\{\frac{2nK_1t^2}{a^2}\Big\}\qquad \textrm{for all}\quad 0\leq t\leq \frac a2.\]
If $\frac a2\leq t\leq t_1$, then
\[\mathbb{E}[e^{tS_n^\ast}]\leq \mathbb{E}[e^{t_1|S_n^\ast|}]\leq  \exp\Big\{n\frac{4\tau_1t_1^q}{a^2}t^2\Big\}.\]
Set $A=\max\{\frac{2K_1}{a^2},\frac{4\tau_1t_1^q}{a^2}\}$. Then
\[\mathbb{E}e^{tS_n}\leq\mathbb{E}e^{tS_n^\ast}\leq e^{nAt^2}\qquad \forall \quad0\leq t\leq t_1.\]
Let
\[\tilde{T_k}=\exp\Big\{tS_k^\ast-kAt^2\Big\}.\]
Then $\{\tilde{T_k},k\geq 1\}$ is a supermartingale. So by Theorem 2.1 in Liu and Watbled \cite{LIU20093101}, we can choose $B>0$ small enough such that for $0\leq x\leq nx_1$,
\begin{equation*}
\begin{aligned}
\mathbb{P}(\max_{1\leq k\leq n}S_k\geq x)&\leq \frac{\mathbb{E}e^{t\alpha \max_{1\leq k\leq n}S_k^\ast}}{e^{t\alpha x}}\\
&\leq e^{-t\alpha x}\mathbb{E}\Big\{\max_{1\leq k\leq n}\tilde{T_k}e^{kAt^2}\Big\}^\alpha\\
&\leq \exp\Big\{-t\alpha x+\alpha nAt^2\Big\}\mathbb{E}\Big\{\max_{1\leq k\leq n}\tilde{T_k}\Big\}^\alpha\\
&\leq (1-\alpha)^{-1}\exp\Big\{-t\alpha x+\alpha nAt^2\Big\}\\
&\leq (1-\alpha)^{-1}\exp\Big\{-\alpha n\sup_{t}\big\{\frac{tx}{n}-At^2\big\}\Big\}\\
&\leq (1-\alpha)^{-1}\exp\Big\{-B\frac{\alpha x^2}{n}\Big\},
\end{aligned}
\end{equation*}
which gives the desired inequality. \hfill\qed
\end{proof*}
\subsection{Bernstain's inequality}
Under Berstein's condition, we obtain the inequalities of Bernstein for NA sequences.
\begin{proposition}
Assume that $\mathbb{E}X_i^2<\infty$. If there exist a constant $M>0$ such that, for any integer $k\geq 2$,
\begin{equation}
\Big|\sum_{i=1}^{n}\mathbb{E}X_i^k\Big|\leq \frac12 k!M^{k-2}B_n,
\end{equation}
where $B_n=\sum_{i=1}^{n}\mathbb{E}X_i^2$. Then, for any $x>0$,
\begin{eqnarray}\label{s2:eq29}
\mathbb{P}(\max_{1\leq k\leq n}S_k\geq x)&\leq& (1-\alpha)^{-1}\exp\biggl\{-\frac{\alpha x^2}{B_n(1+\sqrt{2x\frac{M}{B_n}})+xM}\biggr\}\\
\label{s2:eq30}&\leq& (1-\alpha)^{-1}\exp\biggl\{-\frac{\alpha x^2}{2(B_n+xM)}\biggr\}.
\end{eqnarray}
In particular, we have
\begin{eqnarray}\label{s2:eq100}
\mathbb{P}(\max_{1\leq k\leq n}S_k\geq x)&\leq& 2\exp\biggl\{-\frac{x^2}{2(B_n(1+\sqrt{\frac{2Mx}{B_n}})+Mx)}\biggr\}\\
\label{s2:eq101}&\leq& 2\exp\biggl\{-\frac{x^2}{4(B_n+xM)}\biggr\}.
\end{eqnarray}
\end{proposition}
\begin{proof*}
Let
\[T_k=\exp\big\{tS_k^\ast-\frac{B_kt^2}{2(1-Mt)}\big\}.\]
By Lemma 4.1 in De la Pe\~{n}a \cite{delapea1999}, we have that
\[\mathbb{E}e^{tS_n^\ast}\leq \frac{B_nt^2}{2(1-Mt)}.\]
Hence, $\{T_i,i\geq 1\}$ is a supermartingale. Applying $f(x)=e^{t\alpha x}$ for any $0<t<\frac1M$ and $0<\alpha<1$, to (\ref{s1:eq5}), we have for any $x>0$,
\begin{equation*}
\begin{aligned}
\mathbb{P}(\max_{1\leq k\leq n}S_k\geq x)&\leq \frac{\mathbb{E}e^{t\alpha\max_{1\leq k\leq n}S_k^\ast}}{e^{t\alpha x}}\\
&\leq e^{-t\alpha x}\mathbb{E}\Big\{\max_{1\leq k\leq n}T_ke^{\frac{B_kt^2}{2(1-Mt)}}\Big\}^\alpha\\
&\leq \exp\Big\{-t\alpha x+\frac{\alpha B_nt^2}{2(1-Mt)}\Big\}\mathbb{E}\big\{\max_{1\leq k\leq n}T_k\big\}^{\alpha}\\
&\leq (1-\alpha)^{-1}\exp\Big\{-t\alpha x+\frac{\alpha B_nt^2}{2(1-Mt)}\Big\}.
\end{aligned}
\end{equation*}
Let
\[t=\frac{2Mx^2+x(\sqrt{-B_n(-B_nu^2+2B_nu+2Mx)}+B_nxu)}{B_n^2u+2M^2x^2+B_nMx+2B_nMxu},\]
where $u=(1+\sqrt{\frac{2xM}{B_n}})$. Then we obtain (\ref{s2:eq29}). It is obvious that
\[\sqrt{2xMB_n}\leq \frac{2B_n+Mx}{2},\]
then
\begin{equation*}
\begin{aligned}
\mathbb{P}(\max_{1\leq k\leq n}S_k\geq x)&\leq (1-\alpha)^{-1}\exp\biggl\{-\frac{\alpha x^2}{2B_n+\frac32 Mx}\biggr\}\\
&\leq (1-\alpha)^{-1}\exp\biggl\{-\frac{\alpha x^2}{2(B_n+Mx)}\biggr\},
\end{aligned}
\end{equation*}
which gives (\ref{s2:eq30}). \hfill\qed
\end{proof*}
\subsection{Rio's inequality}
For NA sequences with bounded random variables, we have the following Rio inequality, which is an improved version of the well known Hoeffding-Azuma  inequality.
Recall the following notations of Rio \cite{rio2013}: Let
\[\ell(t)=(t-\ln t-1)+t(e^t-1)^{-1}+\ln (1-e^{-t}) \qquad\text{for all}\quad t>0,\]
and let
\[\ell^\ast(x)=\sup_{t>0}(xt-\ell(t)) \qquad\text{for all}\quad x>0\]
be the Young transform of $\ell(t)$. As quoted by Rio, the following inequality holds
\begin{equation}\label{s2:eq88}
\ell^\ast(x)\geq \max\{2x^2+\frac49 x^4,(x^2-2x)\ln(1-x)\}\qquad \textrm{for}\quad 0\leq x\leq 1.
\end{equation}
\begin{proposition}\label{pr2.9}
Assume that there exist some positive constants $m_i$ and $M_i$ such that
\[
m_i\leq X_i\leq M_i,\qquad\text{for all}\quad 1\leq i\leq n.
\]
Denote
\[M^2(n)=\sum_{i=1}^{n}(M_i-m_i)^2\quad \textrm{and}
\quad D(n)=\sum_{i=1}^{n}(M_i-m_i).\]
Then, for any $t\geq 0$,
\begin{equation}\label{s2:eq89}
\mathbb{E}[e^{tS_n}]\leq \exp\biggl\{\frac{D^2(n)}{M^2(n)}\ell(\frac{M^2(n)}{D(n)}t)\biggr\},
\end{equation}
and, for any $0\leq x\leq D(n)$,
\begin{equation}\label{s2:eq90}
\mathbb{P}(\max_{1\leq k\leq n}S_k\geq x)\leq (1-\alpha)^{-1}\exp\biggl\{-\frac{\alpha D^2(n)}{M^2(n)}\ell^\ast(\frac{x}{D(n)})\biggr\}.
\end{equation}
Consequently, for any $0\leq x\leq D(n)$,
\begin{equation}\label{s2:eq91}
\mathbb{P}(\max_{1\leq k\leq n}S_k\geq x)\leq \biggl(\frac{D(n)-x}{D(n)}\biggr)^{\frac{2D(n)-x}{M^2(n)}\alpha x}.
\end{equation}
In particular, we have for any $0\leq x\leq D(n)$,
\begin{eqnarray}
\mathbb{P}(\max_{1\leq k\leq n}S_k\geq x)&\leq& 2\exp\biggl\{-\frac{D^2(n)}{2M^2(n)}\ell^\ast(\frac{x}{D(n)})\biggr\}\\
&\leq & \biggl(\frac{D(n)-x}{D(n)}\biggr)^{\frac{2D(n)-x}{2M^2(n)}x}.
\end{eqnarray}
\end{proposition}
\begin{remark}
Let us comment on Proposition \ref{pr2.9}.
\begin{itemize}
  \item [1.] It is worth nothing that $\mathbb{E}[X_i]=0$ in Proposition \ref{pr2.9} can be dropped. In fact, since $m_i-\mathbb{E}[X_i]\leq X_i-\mathbb{E}[X_i]\leq M_i-\mathbb{E}[X_i]$ and $M_i-\mathbb{E}[X_i]-(m_i-\mathbb{E}[X_i])=M_i-m_i$, we have, for any $0\leq x\leq D(n)$, the bound (\ref{s2:eq90}) holds for the tail probabilities $\mathbb{P}(\max_{1\leq k\leq n}(S_k-\mathbb{E}[S_k])>x)$.
  \item [2.] Since $(x^2-2x)\ln (1-x)\geq 2x^2$, inequality (\ref{s2:eq89}) implies the following Hoeffding-Azuma inequality
\[\mathbb{P}(S_n>x)\leq \exp\biggl\{-\frac{2x^2}{M^2(n)}\biggr\}.\]
  \item [3.] Taking $\Delta(n)=\max_{1\leq k\leq n}(M_k-m_k)$, we obtain the upper bound: for any $0\leq x\leq n\Delta(n)$,
\[\mathbb{P}(\max_{1\leq k\leq n}S_k>x)\leq \exp\biggl\{-\alpha n\ell^\ast(\frac{x}{n\Delta(n)})\biggr\}\leq \exp\biggl\{-\frac{2\alpha x^2}{n\Delta^2(n)}\biggr\}.\]
\end{itemize}
\end{remark}
\begin{proof*}
Applying $f(x)=e^{tx}$ for any $t\geq 0$, to (\ref{s1:eq2}), we have
\[\mathbb{E}e^{tS_n}\leq \mathbb{E}e^{tS_n^\ast}=\mathbb{E}e^{(tS_n^\ast-t\mathbb{E}S_n^\ast)}.\]
Let $$L(t)=\ell(\Delta_1 t)+\ell(\Delta_2 t)+\cdots+\ell(\Delta_n t),$$ where $\Delta_i=M_i-m_i$. By Lemma 4.1 in Rio \cite{rio2013}, we have
\begin{equation*}
\begin{aligned}
\mathbb{E}e^{tS_n^\ast}&\leq \exp\{L(t)\}\\
&=\exp\Big\{\int_{0}^{t}L'(u)du\Big\}\\
&=\exp\Big\{\int_{0}^{t}\Big(\Delta_1\ell'(\Delta_1u)+\Delta_2\ell'(\Delta_2u)+\cdots+\Delta_n\ell'(\Delta_nu)\Big)du\Big\}.
\end{aligned}
\end{equation*}
By Lemma 4.5 in Rio \cite{rio2013}, $\ell'$ is concave, then we get that
\[\Delta_1\ell'(\Delta_1u)+\Delta_2\ell'(\Delta_2u)+\cdots+\Delta_n\ell'(\Delta_nu)\leq D(n)\ell'\big(\frac{M^2(n)u}{D(n)}\big).\]
So
\begin{equation*}
\begin{aligned}
\mathbb{E}e^{tS_n}&\leq \mathbb{E}e^{tS_n^\ast}\leq \exp\Big\{D(n)\int_{0}^{t}\ell'\big(\frac{M^2(n)u}{D(n)}\big)du\Big\}\\
&=\exp\Big\{\frac{D^2(n)}{M^2(n)}\ell\Big(\frac{M^2(n)t}{D(n)}\Big)\Big\}.
\end{aligned}
\end{equation*}
Let
\[T_k=\exp\big\{tS_k^\ast-\frac{D^2(k)}{M^2(k)}\ell\Big(\frac{M^2(k)t}{D(k)}\Big)\big\}.\]
Then $\{T_k,k\geq 1\}$ is a supermartingale. By (\ref{s1:eq5}), we have for any $0\leq x\leq D(n)$,
\begin{equation*}
\begin{aligned}
\mathbb{P}(\max_{1\leq k\leq n}S_k\geq x)&\leq \frac{\mathbb{E}e^{t\alpha \max_{1\leq k\leq n}S_k^\ast}}{e^{t\alpha x}}\\
&\leq e^{-t\alpha x}\mathbb{E}\Big\{\max_{1\leq k\leq n}T_ke^{\frac{D^2(k)}{M^2(k)}\ell\Big(\frac{M^2(k)t}{D(k)}\Big)}\Big\}^\alpha\\
&\leq \exp\Big\{-t\alpha x+\alpha \frac{D^2(n)}{M^2(n)}\ell\Big(\frac{M^2(n)t}{D(n)}\Big) \Big\}\mathbb{E}\Big\{\max_{1\leq k\leq n}T_k\Big\}^\alpha\\
&\leq (1-\alpha)^{-1}\exp\Big\{-t\alpha x+\alpha \frac{D^2(n)}{M^2(n)}\ell\Big(\frac{M^2(n)t}{D(n)}\Big) \Big\}\\
&\leq (1-\alpha)^{-1}\exp\Big\{-\alpha\sup_{t}\Big(tx-\frac{D^2(n)}{M^2(n)}\ell\Big(\frac{M^2(n)t}{D(n)}\Big)\Big)\Big\}\\
&=(1-\alpha)^{-1}\exp\biggl\{-\frac{\alpha D^2(n)}{M^2(n)}\ell^\ast(\frac{x}{D(n)})\biggr\}.
\end{aligned}
\end{equation*}
By (\ref{s2:eq88}), we get
\begin{equation*}
\begin{aligned}
\mathbb{P}(\max_{1\leq k\leq n}S_k\geq x)&\leq \exp\biggl\{-\frac{\alpha D^2(n)}{M^2(n)}\ell^\ast(\frac{x}{D(n)})\biggr\}\\
&\leq \exp\biggl\{-\frac{\alpha D^2(n)}{M^2(n)}\Big(\big((\frac{x}{D(n)})^2-2(\frac{x}{D(n)})\big)\ln\big(1-(\frac{x}{D(n)})\big)\Big)\biggr\}\\
&=\biggl(\frac{D(n)-x}{D(n)}\biggr)^{\frac{2D(n)-x}{M^2(n)}\alpha x},
\end{aligned}
\end{equation*}
where $0<\alpha<1$. \hfill\qed
\end{proof*}

\end{document}